\documentclass[final,leqno,onefignum,onetabnum]{siamltex1213}
\pdfoutput=1

\usepackage{graphicx}
\usepackage{subfigure}
\usepackage{amsmath}
\usepackage{enumerate}
\usepackage{amssymb}
\usepackage{verbatim}
\usepackage[colorinlistoftodos]{todonotes}
\usepackage{array}

\usepackage{todonotes}
\usepackage{booktabs}
\usepackage{siunitx}
\usepackage{braket}
\usepackage{float}
\usepackage{bm} 
\usepackage{bbm} 
\usepackage[ruled]{algorithm2e}

\newcommand{\edit}[1]{#1}

\graphicspath{ {Figures/} }

\begin{document}


\title{Beyond in-phase and anti-phase coordination in a model of joint action}

\author{Daniele Avitabile 
\thanks{
Centre for Mathematical Medicine and Biology, School of Mathematical Sciences,
University of Nottingham,
University Park, Nottingham, NG7 2RD, UK 
}
\and
Piotr S\l{}owi\'{n}ski 
\thanks{
Department of Mathematics, College of Engineering, Mathematics and Physical
Sciences, \& EPSRC Centre for Predictive Modelling in Healthcare, University of
Exeter, Exeter, Devon, EX4 4QF, UK
}
\and 
Benoit Bardy  
\thanks{
EuroMov, Montpellier University, 700 Avenue du Pic Saint-Loup, 34090 Montpellier,
France \& Institut Universitaire de France, Paris, France
}
\and
Krasimira~Tsaneva-Atanasova 
\thanks{Corresponding author
Email: \texttt{K.Tsaneva-Atanasova@exeter.ac.uk}}
Department of Mathematics, College of Engineering, Mathematics and Physical
Sciences, \& EPSRC Centre for Predictive Modelling in Healthcare, University of
Exeter, Exeter, Devon, EX4 4QF, UK
}

\graphicspath{{./Figures/}}

\maketitle

\begin{abstract}
In 1985 Haken, Kelso and Bunz proposed a system of coupled nonlinear oscillators as a
model of rhythmic movement patterns in human bimanual coordination. 
Since then, the Haken-Kelso-Bunz (HKB) model has become a modelling paradigm applied
extensively in all areas of movement science, including interpersonal motor
coordination.  However all previous studies \edit{have followed a line of analysis
based on slowly varying amplitudes and rotating wave approximations. These
approximations lead to a reduced system, consisting of} a single differential
equation representing the evolution of the relative phase of the two coupled
oscillators: the HKB model of the relative phase Here we take a different approach
and systematically investigate the behaviour of the HKB model in the full
four-dimensional state space \edit{and for general coupling strengths}. We perform
detailed numerical bifurcation analyses and reveal that the HKB model supports
previously unreported dynamical regimes as well as bi-stability between a variety of
coordination patterns.  Furthermore we identify the stability boundaries of distinct
coordination regimes in the model and discuss the applicability of our findings to
interpersonal coordination and other joint-action tasks.
\end{abstract}


\section{Introduction}

Many body movements are periodic in their nature \cite{Fuchs2008}. For example postural
sway \cite{Bardy2002}, walking \cite{Collins1993,Golubitsky1999}, running \cite{Laffaye2007}, swimming \cite{Seifert2004}, galloping \cite{Collins1993,Golubitsky1999} and juggling
\cite{Zelic2012} have a cyclic pattern in the position of the end effectors or joint angles.
Synchronisation is a fundamental aspect of oscillatory coordination
dynamics in human and animal body movements \cite{Kelso1997} and has been found in
many different situations \cite{Saltzman1987}.
Coordination is characterised by a bounded temporal relationship created
by a convergent dynamical process \cite{Jirsa2005,Mortl2012}. 
Coordination regimes depend on \edit{symmetries} and couplings between oscillators.
Frequency entrainment, where two oscillators adopt a central frequency, occurs
even with a very weak coupling. \edit{With a relatively strong coupling or if the
system is symmetric}, phase entrainment can also take place. These processes
may be continuous or intermittent, that is the phases of the two oscillators
may also align periodically \cite{Fuchs2008,Pikovsky2002,Kurths2003}.

In the case of two coupled oscillators the regular patterns of coordination are well captured by the properties of
the relative phase between the periodic movements of the two subsystems \cite{Kelso1987,Kelso1997}.
The simplest pattern is observed when the phase of the two oscillators coincide
to give in-phase ($0^{\circ}$) monostable coordination pattern. An example \edit{of
this behaviour is given by} iso-lateral limb movements \cite{Buchanan2006}.
Monostable anti-phase
($180^{\circ}$)
coordination can also occur and an example of such behaviour is observed
in team sports (competitive games) \cite{Bourbousson2010,Esteves2012,Duarte2012,Davids2013}. In many real
systems, anti-phase stability coexists with in-phase stability
\cite{Haken1985,Kelso1997,Mortl2012,Warren2006}. 
\edit{Previous studies address the modelling the two coupled oscillators as a nonlinear dynamical
system, the fitting of its periodic orbits to human movements~\cite{Kay1987}, and
the systematic analysis of the effects of linear and nonlinear terms to the observed
limit cycles~\cite{Beek1995}.}
The observed relations between frequency and amplitude
\cite{Haken1985} as wells as peak velocities \cite{Mortl2012} in many but not all
\cite{Beek2002} oscillatory movements turned out to be well represented by a hybrid
oscillator \cite{Haken1985} formed by a combination of Van der Pol and Rayleigh
nonlinear damping
terms.

A classical example of model exhibiting bi-stability is the so-called \edit{HKB model proposed in
the seminal work by Haken, Kelso and Bunz \cite{Haken1985,Fuchs2008}. The model, which was originally developed for bimanual
finger coordination \cite{Kelso1984}, has found to be representative of a wide range of applications in
human movement \cite{Fuchs2008,Calvin2011} suggesting that the dynamics it produces
are somehow fundamental and make formal construct for the study of coordination
dynamics \cite{Kelso1987,Kelso1997}. Although the model was originally developed in
order to account for an intra-personal phenomena, the same patterns have been shown to
be representative of both sensorimotor and 
interpersonal behaviours \cite{Kelso1990,Schmidt1990,Schmidt2008}.} The model successfully
reproduces not only the patterns of stability observed in bimanual coordination
experiments but also their dependence upon frequency \cite{Haken1985}. 
The HKB model admits a potential function that yields the
experimentally observed change in attractors' landscapes. 
\edit{Furthermore, the HKB model and its stochastic extension reproduced the
characteristic fluctuation increase and slowing down observed experimentally near
instabilities~\cite{Schoener1986}.}

The development of the HKB model has been inspired by the in-phase and anti-phase
coordination dynamics observed in bimanual coordination in the context of the finger
movements experiment \cite{Kelso1981,Kelso1983},\cite{Kelso1984}. Therefore most previous research has
focused on a fixed set of model parameters that guarantees the stability of these particular dynamics.
Furthermore, significant contributions to understanding these coordination patterns
(\edit{albeit in a narrow parameter range and with limiting assumptions on the
parameters controlling the coupling strength}) have been made for
different oscillator frequencies and inputs
\cite{Fuchs1996,Fuchs2000a,Jirsa2000,Beek2002,Assisi2005,Calvin2011} as well as noise
in the system \cite{Schmidt1990,Daffertshofer1999,Schmidt2008}.  \edit{All
previous mathematical analyses of the HKB model have focused on the relative phase
dynamics, under the assumption that the amplitude of the coupled oscillators is constant
\cite{Haken1985,Fuchs1996,Daffertshofer1999,Fuchs2000a,Beek2002,Assisi2005}}.  Several
recent articles have studied the phase-approximation dynamics in the HKB model by
considering the multiple stable states of the system and the ability to switch
between them by changing the frequency and the coupling parameters 
\cite{Leise2007}. The bifurcations leading to transitions between anti-phase and
in-phase dynamics in a reduced phase approximation of the HKB model \cite{Frank2012}
have been also studied. To our knowledge, however, a bifurcation analysis of the full
four-dimensional HKB model, considering all model parameters \edit{as well as general
(i.e. weak and strong) coupling strengths}, has not been performed.
Such analysis could provide an insight into other possible qualitative behaviours that
the solutions of the model might exhibit, as well as characterise the possible
changes in the dynamics of the solutions corresponding to any changes in the
parameter values of the model. We also note recent further developments of dynamical
systems' approaches for studying sensorimotor dynamics, involving dynamical
repertoires, hierarchies of timescales and structured flows on
manifolds~\cite{Huys2014}.

\edit{Given that the HKB model is a widely accepted tool in this field, it is
imperative to examine systematically all the possible coordination regimes supported
by this system. In addition, classifying changes of dynamical regimes in terms of
positions and velocities of the two coupled oscillators would undoubtedly shed light
on the HKB model's applicability to explain movement coordination in joint actions
and human interactions with an adaptive virtual partner
(VP)\cite{Kelso2009,Mattout2012,Dumas2014,ZhaiSMC2014,ZhaiCDC2014,Zhai2016}.}
\edit{In the present paper} we take a different approach in analysing the HKB model,
as we study the full four
dimensional system of first order differential equations describing the evolution of
the positions and velocities of the two coupled oscillators.  We begin by
characterising the local and global dynamics of the single HKB oscillator and reveal a
global transition in the model that governs the existence of periodic solutions in a
range of \edit{the oscillator's} parameter values. We proceed by systematically
characterising the full HKB model dynamics not only by varying the coupling strength
parameters but also the rest of the model parameters, i.e. the parameters governing
the single oscillator's properties.  In addition to the very well studied
coordination patterns we find a stable phase-locked solution that spans a wide range
of relative phases and persists for a wide range of model parameters' values.  We also
show that relaxing the constant amplitude assumption allows for much richer
co-ordination dynamics and co-existence of various stable coordination attractors
(multi-stability regimes).

\section{Results}

\subsection{Intrinsic properties of the oscillator in the HKB model}

Recently, a significant scientific effort has been put towards the development of VP interaction systems. In particular the single HKB oscillator is being used to drive the movement dynamics exhibited by the VP \cite{Kelso2009,Dumas2014,ZhaiSMC2014,ZhaiCDC2014,Zhai2016}. \edit{The dynamics of the model is an important consideration in designing such systems and in particular for parametrising the ordinary differential equation that governs the behaviour of the VP. For example depending on the constraints of the experimental set-up, a certain range of amplitude and/or frequency for the VP periodic behaviour might be desirable. Although, some properties of the HKB oscillator have been measured and studied both experimentally and analytically \cite{Kay1987, Kay1991}, the dynamics of the
single HKB oscillator has not been systematically investigated theoretically.
To address this gap we begin by examining a single HKB oscillator:}

\begin{equation}
\begin{aligned}
  \ddot{x} &= - \dot{x} \left( \alpha x^2  + \beta  \dot{x}^2 - \gamma \right)-\omega ^2 x, \nonumber
\end{aligned}
\end{equation}
which could be written as a planar autonomous dynamical system of the form:
\begin{equation}
\begin{aligned}
  \dot{x} &= y, \\
  \dot{y} &= - y \left( \alpha x^2  + \beta  y^2 - \gamma \right)-\omega ^2 x,
\end{aligned}
\label{eq:HKB}
\end{equation}
where $x$ represents the position, $y$ the velocity, $\omega \in \mathbb{R}^+$ is related to the natural frequency of the oscillator
and $\alpha, \beta, \gamma \in \mathbb{R}$ are parameters governing the intrinsic dynamics of equation \eqref{eq:HKB}.
\begin{figure}
  \centering
  \includegraphics{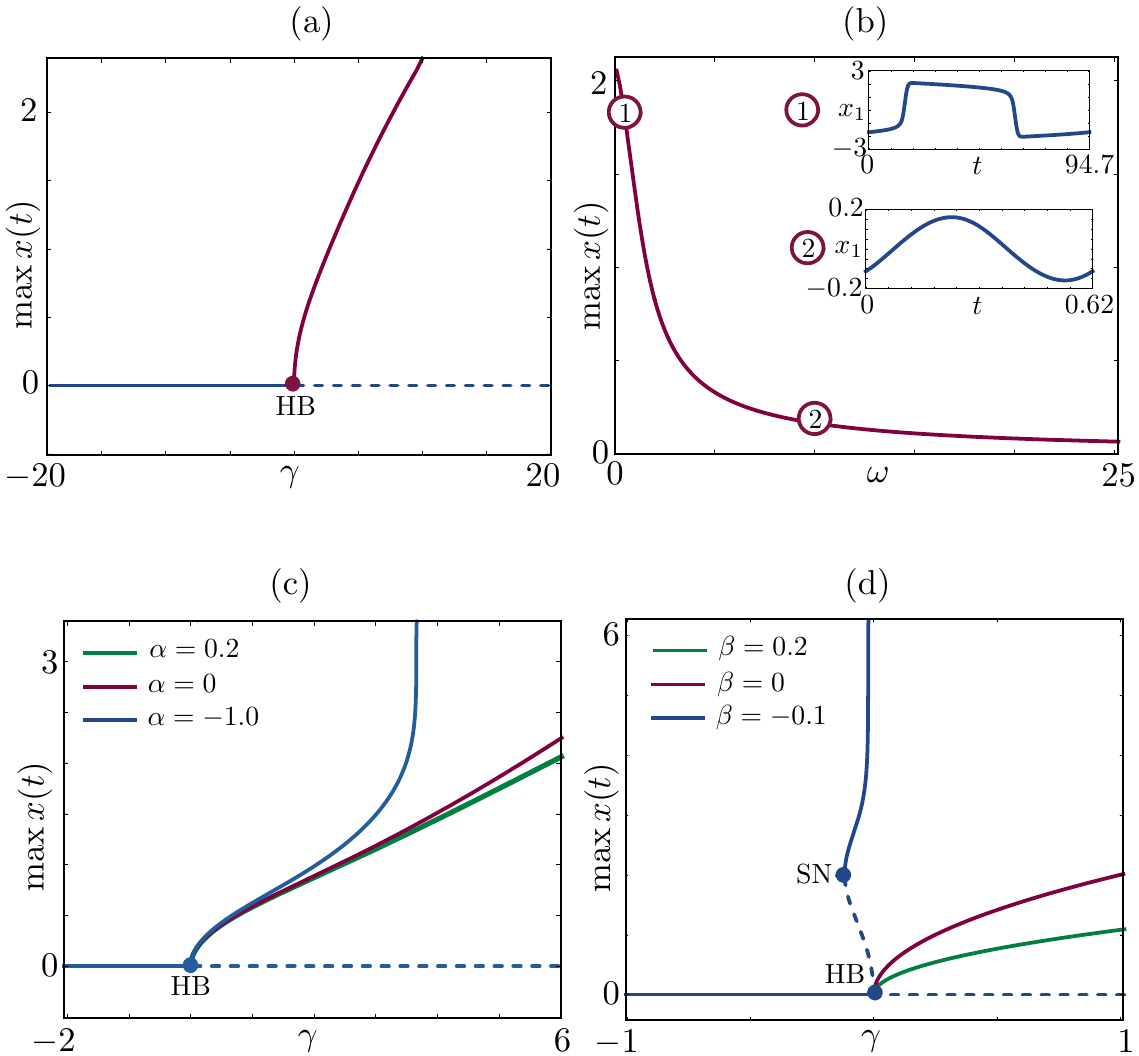}
  \caption{Bifurcation diagrams for a single HKB oscillator. (a): the trivial
  equilibrium becomes unstable at a supercritical Hopf bifurcation (HB) in the
  continuation parameter $\gamma$ for $\omega=2$, $\alpha=1$, $\beta=1$. (b): the
  periodic orbit for $\gamma = 2$ is continued in the parameter $\omega$. The lower
  $\omega$, the larger the oscillations amplitude and the longer the period.
  (c)--(d): continuations in $\gamma$ are repeated for various values of $\alpha$ and
  $\beta$. In panel (c), for $\alpha=-1$, the periodic branch undergoes a global bifurcation
  (vertical asymptote), whereas in panel (d), for $\beta = -0.1$, the Hopf bifurcation is
  subcritical, and the emanating branch restabilises at a saddle-node
  bifurcation, before disappearing in a global bifurcation. 
  \edit{Solid (dashed) lines represent stable (unstable) states of \eqref{eq:HKB}.}}
  \label{fig:singleOscillator}
\end{figure}

\edit{The single HKB oscillator is a hybrid Rayleigh~-~Van der Pol \cite{Haken1985} planar system and although the analysis of planar systems of ordinary differential equations is very well established \cite{Mook1979,Jordan2007,Guckenheimer2013}, it has not been applied to the single HKB oscillator model. Furthermore whenever planar systems are coupled, they are often studied in the weak coupling limit, which we don't require for the numerical-continuation analysis presented here.}
In our analysis, we focus on the global dynamics of the system and aim to
characterise all possible dynamic states that the single HKB oscillator model
supports, as well as their dependence on all model parameters. System \eqref{eq:HKB}
admits the origin $(0,0)$ as a trivial steady state for any parameter value $\omega \in \mathbb{R}^+$.
Given $\omega>0$, the Jacobian matrix at the trivial equilibrium $(x,y)=(0,0)$ is
\begin{equation*}
J = 
\begin{bmatrix}
0 & 1 \\
-\omega^2 & \gamma
\end{bmatrix}
\end{equation*}
For $|\gamma| \geq 2\omega$, the Jacobian has a pair of non-zero real eigenvalues:
\begin{equation}
\begin{aligned}
\lambda = \frac{\gamma \pm \sqrt{\gamma^2 - 4\omega^2}}{2} \nonumber
\end{aligned}
\end{equation}
Thus, the equilibrium is a stable node (sink) for $\gamma<0$ and unstable node (source) for $\gamma>0$.
For $|\gamma| < 2\omega$, the Jacobian has a pair of complex conjugate eigenvalues of the form:
\begin{equation}
\begin{aligned}
\lambda = \frac{\gamma}{2} \pm i\frac{\sqrt{4\omega^2-\gamma^2}}{2} \nonumber
\end{aligned}
\end{equation}
\edit{Hence, the equilibrium is a stable focus (spiral sink) for $-2\omega<\gamma<0$ and unstable focus (spiral source) for $0<\gamma<2\omega$.}

Changing the value of the parameter $\gamma$ near $\gamma=0$ leads to a change
in the sign of the eigenvalues' real part, which is associated with loss or gain
of stability. The system undergoes a Hopf bifurcation at $\gamma=0$, which gives
rise to oscillations. 
We could analytically verify further the sufficient conditions for the existence
of Hopf bifurcation by showing that:
\begin{align*}
  \frac{\partial \lambda_\textrm{r}(\gamma)}{\partial \gamma}\Bigg \vert_{\gamma=0} & = \frac{1}{2} \neq 0, \\
  l_1(\gamma) \vert_{\gamma=0} &= 
  \frac{-(\alpha+3\beta \omega^2 )}{2\omega(\omega^2+1)} \neq 0 \quad \iff \quad \alpha+3\beta \omega^2 \neq 0,
\end{align*}
where $\lambda_\textrm{r}$ and $l_1$ are the real part of the
eigenvalues and the first Lyapunov coefficient \cite{Kuznetsov2013}, respectively. 
\edit{The sign of the first Lyapunov coefficient \cite{Kuznetsov2013} determines
whether the Hopf bifurcation is subcritical or supercritical, hence we are in the 
supercritical (subcritical) case if $\alpha+3 \beta \omega^2>0$ ($<0$).}
The system \eqref{eq:HKB} has a degenerate Hopf bifurcation when $\alpha+3 \beta
\omega^2=0$.

Next we carry out bifurcation analysis using numerical continuation in AUTO
\cite{AUTO_orig}. We set $\texttt{NTST}=50$, $\texttt{NCOL}=4$ for the mesh,
and $\texttt{EPSL}=10^{-9}$, $\texttt{EPSU}=10^{-9}$ for the tolerances of the Newton solver.
In Figure~\ref{fig:singleOscillator} (a) we continue the trivial steady state $(x,y)=(0,0)$ in $\gamma$: oscillations
arise at a supercritical Hopf bifurcation at $\gamma = 0$ since for $\alpha=1$
and $\beta=1$, $\alpha+3 \beta \omega^2>0$. Stable periodic solutions exist for
various values of the intrinsic frequency $\omega$: the lower $\omega$, the
larger the oscillations amplitude and the longer the period (see inset in
Figure~\ref{fig:singleOscillator} (b)). Similar scenarios are found for various
combinations of $\alpha$ and $\beta$ (Figure~\ref{fig:singleOscillator}c--d). 
The amplitude of the periodic solutions increase as either $\alpha$ or $\beta$ are decreased.
For $\alpha=-1$ the Hopf bifurcation is supercritical and the oscillatory branch is stable,
whereas for $\beta=-0.1$ the Hopf bifurcation is subcritical and oscillatory
branches, which are originally unstable, re-stabilise at a saddle-node. We note
that system \eqref{eq:HKB} exhibits bi-stability between a stable equilibrium and a stable periodic states
in the case of $\beta=-0.1$ (Figure~\ref{fig:singleOscillator}d).

The above analysis reveals that when parameters $\alpha$ and $\beta$ have opposite signs
(Figure~\ref{fig:singleOscillator} (c) and (d)) there is a critical value for
$\gamma$ at which the amplitude (and period) of the stable limit cycle solutions
in the model rapidly increase to infinity. As a result all periodic solutions vanish
for values of $\gamma$ above this critical value. Furthermore, such transitions occur robustly for a large range 
of $\alpha$ and $\beta$ parameter values. \edit{We believe it is important to
understand where and why this singularity occurs, as it corresponds to a non-physical
behaviour. Since this feature has not previously been reported in the literature on
the HKB model, we present a thorough investigation of this phenomenon.}
In order to analyse the behaviour of the system at infinity we employ methods
presented in Chapter 3.10 of reference \cite{Perko}. We start by projecting system
\eqref{eq:HKB} on the Poincar\'e sphere using the following transformation:
\begin{eqnarray*}
&X=\dfrac{x}{\sqrt{1+x^2+y^2}}, \;\; Y=\dfrac{y}{\sqrt{1+x^2+y^2}}, \;\; Z=\dfrac{1}{\sqrt{1+x^2+y^2}},
\end{eqnarray*}
which defines one-to-one correspondence between points $(X,Y,Z)$ on the upper hemisphere $\mathrm{S}^2$ with $Z>0$ and points $(x,y)$ in the plane defined by:
\begin{eqnarray*}
&x=\dfrac{X}{Z}, \;\; y=\dfrac{Y}{Z},\\
\end{eqnarray*}
The points on the equator of $\mathrm{S}^2$ correspond to points at infinity of $\mathbb{R}^2$. 
Under the transformation above the HKB oscillator on $\mathrm{S}^2$ with $Z>0$ is given by:
\begin{equation}
\begin{aligned}
&\dot X  =\dfrac{Y}{Z^2} \big[ \alpha X^3 Y+\beta X Y^3+Z^2 \big(-\gamma X Y+(\omega^2-1) X^2+1\big) \big]\\
&\begin{split}
  \dot Y   = \dfrac{1}{Z^2} \Big[ & (Y^2-1) Y(\alpha X^2 +\beta Y^2) \\
                                  & + Z^2 \big(Y (-\gamma (Y^2-1)-X Y)+\omega^2 X (Y^2-1)\big)
	                    \Big] 
\end{split} \\
&\dot Z =\dfrac{Y}{Z} \big[\alpha X^2 Y+\beta Y^3+Z^2 \big((\omega^2-1) X-\gamma Y\big)\big]
\end{aligned}
\label{eq:sphere_HKB}
\end{equation}

System \eqref{eq:HKB} has 8 equilibria on the equator $X^2+Y^2=1$ of $\mathrm{S}^2$ (see Theorem 1 from Chapter 3.10 of \cite{Perko}) that represents the limit $x, y \rightarrow \infty$. In general, the equilibria are given by the solutions of the following equation:
\begin{equation}
XQ_m(X,Y)-YP_m(X,Y)=0
\label{eq:theorem1}
\end{equation}
where $P_m$ and $Q_m$ are homogeneous $m$-th degree polynomials in $x$ and $y$  
according to the following representation of the system \eqref{eq:HKB}:
\begin{eqnarray}
\begin{aligned}
\dot x&=P(x,y)=P_1(x,y)+\dots+P_m(x,y)\\
\dot y&=Q(x,y)=Q_1(x,y)+\dots+Q_m(x,y)
\end{aligned}
\end{eqnarray}
In our case the highest degree homogeneous polynomials are:
\begin{eqnarray}
\begin{aligned}
P_3(x,y)&=0\\
Q_3(x,y)&=-\alpha x^2y-\beta y^3 
\end{aligned}
\end{eqnarray}
Hence, all equilibria at the equator are the solutions of the following system of equations:
\begin{eqnarray}
\begin{aligned}
&X^2+Y^2=1\\
&XQ_3(X,Y)-YP_3(X,Y)=-\alpha X^3Y-\beta Y^4=0
\end{aligned}
\label{eq:theorem1_case}
\end{eqnarray}
and are given by:
\begin{eqnarray}
\begin{aligned}
X&=0, &Y&=\pm1,\\
X&=\pm1, &Y&=0,\\
X&=\pm \dfrac{\sqrt{\alpha}}{\sqrt{\alpha-\beta}}, &Y&=\pm \dfrac{\sqrt{\beta}}{\sqrt{\beta-\alpha}.}
\end{aligned}
\label{eq:equilibria}
\end{eqnarray}
The flow between the nodes is determined using the following equation (see Theorem 1 from Chapter 3.10 of \cite{Perko}):
\begin{equation}
G_{m+1}=\cos\theta Q_m(\cos\theta,\sin\theta)-\sin\theta P_m(\cos\theta,\sin\theta)=0,
\end{equation}
where $\theta$ is an angle along the equator.

\begin{figure}
\centering
\includegraphics[width=\textwidth]{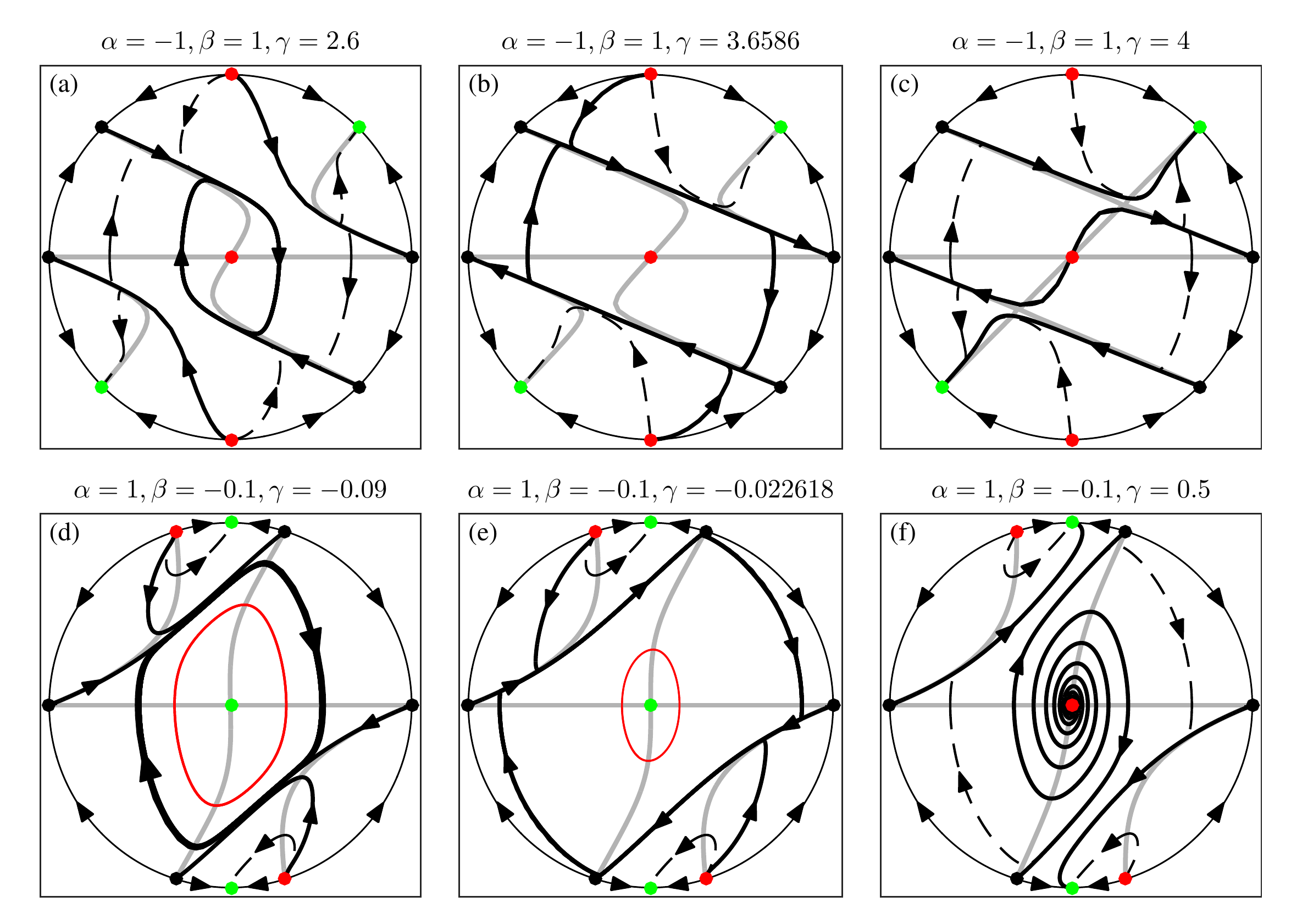}
\caption{Global phase portraits of the system Eq.~(\ref{eq:sphere_HKB}), projected on the $(X,Y)$-plane, for different parameters values. Green dots indicate stable equilibria; red dots indicate unstable equilibria; black dots indicate equilibria of a saddle type; red line indicates unstable periodic orbit; thick black lines indicate heteroclinic connections between different equilibria or between equilibria and stable periodic orbits; grey lines indicate nullclines; dashed lines examples of trajectories; arrows indicate direction of the flow.}
\label{fig:disks}
\end{figure}
 
The flow between the equilbiria on the equator of the Poincar\'e sphere is counter-clockwise if $G_{m+1}>0$ and clockwise where $G_{m+1}<0$. 
We find that only the equilibria $X=0, Y=\pm1$ are hyperbolic. They are stable
nodes for $\beta<0$ and are un-stable nodes for $\beta>0$. The other six
equilibria as given in \eqref{eq:equilibria} are non-hyperbolic. We established
their types by combining information gathered from the flow on the equator and
from numerical integration of the transformed system \eqref{eq:sphere_HKB}.  
We summarise our findings in two representative cases in which, as the parameter
$\gamma$ increases, the period and amplitude of the stable periodic orbit grows
to infinity exponentially fast and the periodic orbit disappears. More
specifically, at the critical value $\gamma^*$, the stable periodic orbit
becomes a heteroclinic cycle connecting four equilibria of saddle type at the
equator on the Poincar\'e sphere.

In Figure \ref{fig:disks} we illustrate how the structure of the global phase
portrait of the system \eqref{eq:sphere_HKB}, projected on the $(X,Y)$-plane,
changes with increasing $\gamma$. 
Figures \ref{fig:disks}(a)--(c), for $\alpha=-1, \beta=1$, show the transition occurring as $\gamma$ 
is varied in the bifurcation diagram of Fig.~\ref{fig:singleOscillator}c,  
Figures \ref{fig:disks}(d)--(f), for $\alpha=1, \beta=-0.1$, show the transition occurring as $\gamma$ 
is varied in the bifurcation diagram of Fig.~\ref{fig:singleOscillator}d.
In both cases the disappearance of the stable limit cycle solution in the model is
due to the same mechanism. However, depending on the signs of the parameters
$\alpha, \beta$, different invariant objects are involved in the transition.
Panels (a)-(c) in Fig. \ref{fig:disks}
show that there are two types of connecting orbits in the
phase space of the HKB oscillator. The first type connects the unstable
equilibria $(0,\pm1)$ (red dots) with the saddle points $(\pm1,0)$ (black dots)
and the second connects the saddle points $(\mp
\sqrt{\alpha}/\sqrt{\alpha-\beta}, \pm \sqrt{\beta}/\sqrt{\beta-\alpha})$ (black
dots) with the stable periodic orbit surrounding the unstable equilibrium at the
origin $(0,0)$. As the parameter $\gamma$ increases, the two types of connections
become tangent and the periodic orbit stretches along the $X$ axis as depicted
in panel (b) for value of $\gamma=3.65860608978$ (just before the transition).
At the critical value, $\gamma=\gamma^*$, the periodic orbit becomes a
heteroclinic cycle connecting four saddle equilibria. After the transition the
heteroclinic cycle disappears and the global phase portrait changes. In panel
(c) we show that after the transition there are connections between the saddle
points $(\mp \sqrt{\alpha}/\sqrt{\alpha-\beta}, \pm
\sqrt{\beta}/\sqrt{\beta-\alpha})$ and the stable equilibria $(\pm
\sqrt{\alpha}/\sqrt{\alpha-\beta}, \pm \sqrt{\beta}/\sqrt{\beta-\alpha})$, and
between the saddle points $(\pm1,0)$ and the unstable equilibrium at the origin
$(0,0)$. In this case the single HKB oscillator has stable periodic solutions
only for $\gamma \in (0, \gamma^*)$.
Panels (d)-(e) in Fig. \ref{fig:disks} demonstrate that, for $\alpha=1,
\beta=-0.1$, in addition to the stable periodic orbit there is also an unstable
periodic orbit (red loop) surrounding the stable equilibrium at the origin
$(0,0)$ (green dot). \edit{Although unstable periodic orbits could not be observed experimentally, such objects are important from dynamical systems point of view.
For example in this case the branch of unstable periodic orbits forms the boundary between the basins of attraction of the coexisting stable equilibrium and stable periodic orbit for $\gamma \in (\gamma_\textrm{SN}, \gamma^*)$ (see panel (d) in Fig. \ref{fig:singleOscillator}.}
Irrespective of the presence of unstable limit cycle we
find again two types of connecting orbits in the phase space for
$\gamma<\gamma^*$ as shown in panel (d). The first type connects the unstable
equilibria $(\mp \sqrt{\alpha}/\sqrt{\alpha-\beta}, \pm
\sqrt{\beta}/\sqrt{\beta-\alpha})$ (red dots) with the saddle points $(\pm
\sqrt{\alpha}/\sqrt{\alpha-\beta}, \pm \sqrt{\beta}/\sqrt{\beta-\alpha})$ (black
dots). The second type connects the saddle points $(\pm1,0)$ (black dots) with
the stable periodic orbit. 
Here we observe again that the connections become tangent to the periodic
orbit, as it stretches along the $X$-axis growing into a heteroclinic cycle between
four saddle equilibria (black dots) for $\gamma=\gamma^*$, as depicted in panel
(e) where $\gamma=-0.022618$ (just before the transition). After the transition
$\gamma>\gamma^*$  the heteroclinic cycle disappears and the invariant objects
of the system reconnect. This, however, occurs in a different manner compared to
the case presented in panels (a)-(c). The saddle equilibria $(\pm1,0)$ are now connected
with stable nodes $(0,\pm1)$ and the unstable periodic orbit is connected to the
saddle points $(\pm \sqrt{\alpha}/\sqrt{\alpha-\beta}, \pm
\sqrt{\beta}/\sqrt{\beta-\alpha})$. In panel (f) we show the phase portrait for
$\gamma=0.5$, which illustrates the connections after the unstable periodic
orbit disappeared in a subcritical Hopf bifurcation (at $\gamma=0$, compare with
Fig.~\ref{fig:singleOscillator}). In this case the single HKB oscillator has
stable periodic solutions only for $\gamma \in (\gamma_{SN},\gamma^*)$.

\subsection{Bifurcation analysis of the full HKB model}

\subsubsection{Full system model equations}

Previous analysis of the HKB model has focussed on the dynamics of the relative phase
that is given by the difference of the two oscillators' phases.
However, in applications involving VP interaction environments
\cite{Dumas2014,ZhaiSMC2014,ZhaiCDC2014},
other properties of the HKB
model dynamics become crucial. Such properties include the amplitude and phase of the
oscillatory solutions, as well as their existence, parameter dependence and
stability. In order to address these questions we focus below on the full HKB system.
The original HKB model evolves in time (measured in seconds) according to a set of nonlinear differential equations \cite{Haken1985}:

\begin{eqnarray}
\label{eq1}
\ddot{x_1} + \dot{x_1} \left( \alpha x_1^2  + \beta  \dot{x_1}^2 - \gamma \right)+\omega ^2 x_1 &=& I_{12}(\dot{x_1}, \dot{x_2}, x_1, x_2)\\
\nonumber
\ddot{x_2} + \dot{x_2} \left( \alpha x_2^2  + \beta  \dot{x_2}^2 - \gamma \right)+\omega ^2 x_2 &=& I_{21}(\dot{x_1}, \dot{x_2}, x_1, x_2),
\end{eqnarray}
where $x_1$ and $x_2$ represent the position of the two agents' end effectors and
\begin{eqnarray}
\label{eq2}
I_{12}(\dot{x_1}, \dot{x_2}, x_1, x_2) &=& (a+b(x_1-x_2)^2)(\dot{x_1}-\dot{x_2})\\
\nonumber
I_{21}(\dot{x_1}, \dot{x_2}, x_1, x_2) &=& (a+b(x_2-x_1)^2)(\dot{x_2}-\dot{x_1}),
\end{eqnarray}
are coupling functions with coefficients $a, b \in \mathbb{R}$.
The above system of two coupled second order ordinary differential equations (ODEs) \eqref{eq1} can be written as a four dimensional autonomous system of first order ODEs: 
\begin{align}
\label{eqn4}
  \dot{x_1} &= y_1 \\
  \nonumber
  \dot{x_2} & = y_2 \\
  \nonumber
  \dot{y_1} & = (a+b(x_1-x_2)^2)(y_1-y_2) - (y_1 \left( \alpha x_1^2  + \beta  {y_1}^2 - \gamma \right)+\omega ^2 x_1)\\
  \nonumber
  \dot{y_2} & = (a+b(x_2-x_1)^2)(y_2-y_1) - (y_2 \left( \alpha x_2^2  + \beta  {y_2}^2 - \gamma \right)+\omega ^2 x_2 ),
\end{align}
\edit{where $x_i$ and $y_i$ represent position and velocity of the $i$th agent's end
effector, respectively.
The resulting dynamical system has a four-dimensional state space \cite{Calvin2011}. 
The parameter $\omega$ (commonly referred to as {\it eigenfrequency}) defines, in 
conjunction with $\alpha, \beta$ and $\gamma$, the intrinsic dynamics of the two
coupled oscillators. The oscillators' positions and velocities are coupled via the
parameters $a$ and $b$, commonly referred to as {\it coupling strengths}.} 
The HKB model behaviour then depends on the intrinsic dynamics parameters as well as
the coupling strengths. Although coordination/synchronisation in system \eqref{eqn4}
emerges as a consequence of coupling, its dynamics (i.e., number, type and stability
of coordination patterns) depends not only on the nature of the coupling but also on
the intrinsic properties of each coupled oscillator. \edit{In the HKB system
\eqref{eqn4}, both the intrinsic dynamics and the couplings are highly nonlinear,
opening up the possibility of obtaining multistability and hence multifunctionality.}

\subsubsection{Coordination regimes in the HKB model}

\begin{figure}
  \centering
  \includegraphics{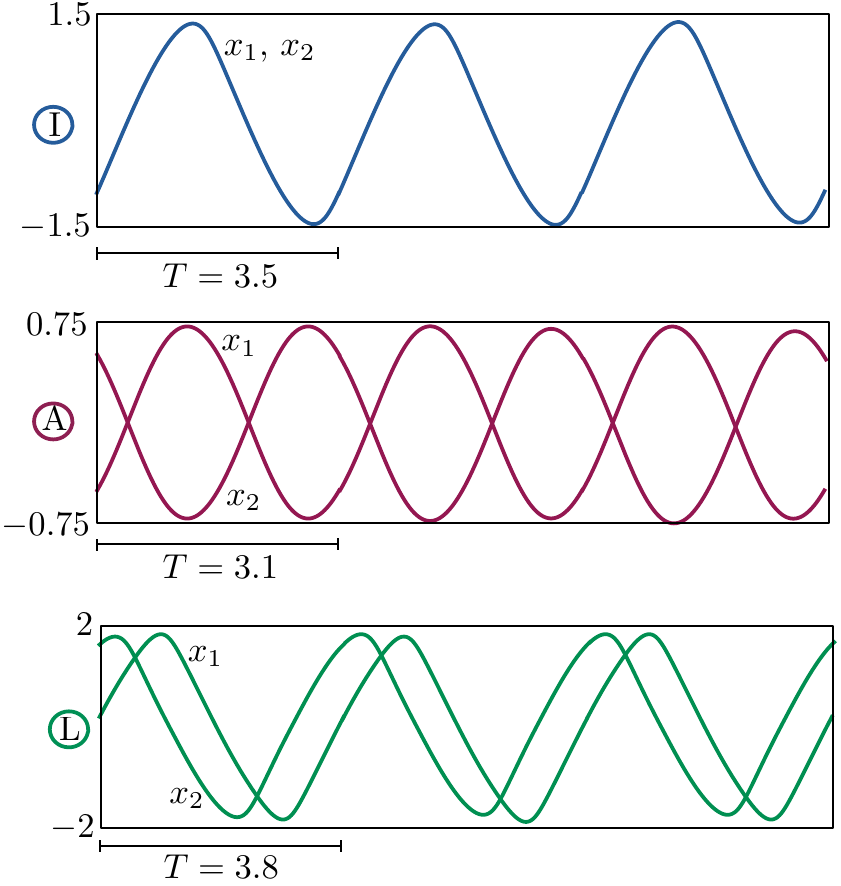}
  \caption{Examples of stable in-phase (I), anti-phase (A) and phase-locked (L) solutions.
  Solutions and parameter values are also indicated in the bifurcation diagrams of
  Figure~\ref{fig:gammaContinuationCompound}.}
  \label{fig:solutionProfiles}
\end{figure}

\begin{figure}
  \centering
  \includegraphics[width=\textwidth]{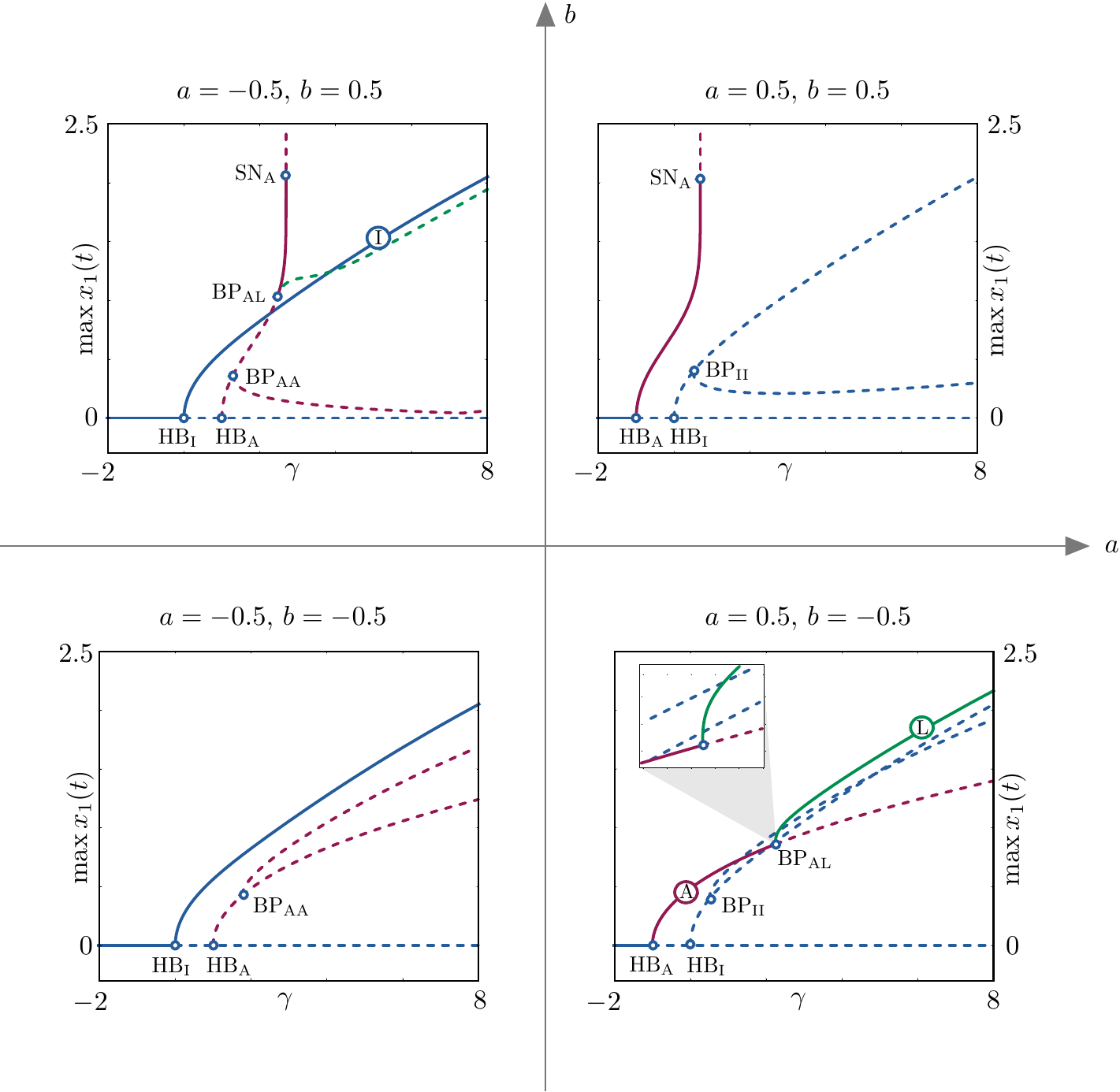}
  \caption{Representative bifurcation diagrams in \edit{the parameter} $\gamma$ for
  all possible combinations of coupling strengths, $a$ and $b$. \edit{Solid (dashed)
  lines represent stable (unstable) states of \eqref{eqn4}.}}
  \label{fig:gammaContinuationCompound}
\end{figure}

In this section we study the existence and stability of the possible coordination regimes in the full HKB model~\eqref{eqn4} by conducting a
systematic analysis in all model (control) parameters. 
The numerical bifurcation analysis
is carried out using numerical continuation in AUTO \cite{AUTO_orig}. 
We set $\texttt{NTST}=50$, $\texttt{NCOL}=4$ for the mesh,
and $\texttt{EPSL}=10^{-9}$, $\texttt{EPSU}=10^{-9}$ for the tolerances of the Newton solver.
We perform time-stepping simulations of the model \eqref{eqn4} in
\textsc{MATLAB} \cite{MATLAB:2012}, using the \texttt{ode45} solver with default numerical settings.
In the simulations presented below we use the following typical intrinsic dynamics parameter values as default, $\alpha=1,  \beta=1, \gamma=1$ and $\omega=0.2$, unless otherwise stated in the figure legends.

In agreement with previously performed analysis on the HKB relative phase dynamics \cite{Haken1985,Fuchs1996,Daffertshofer1999,Fuchs2000a,Beek2002,Assisi2005}
we confirm existence and study the stability of the well characterised in-phase and anti-phase oscillatory solutions.
Moreover we find a new family of stable periodic phase-locked solutions characterised
by relative phase in the interval $(0^{\circ}, 180^{\circ})$. \edit{These solutions
are found to be stable in a wide range of parameter values}.
We note that this family of solutions is unstable for the commonly used set of model parameters based on \cite{Haken1985}. 
Examples of the three solution types described above are plotted in Figure~\ref{fig:solutionProfiles}.
We show how such solutions are born when we vary $\gamma$ 
for various combinations of the parameters $a$ and
$b$ in the bifurcation diagrams of Figure~\ref{fig:gammaContinuationCompound}, whose
branches are color coded as in
Figure~\ref{fig:solutionProfiles}. Here and henceforth, we use subscripts I, A, L (or
combinations thereof) to indicate bifurcations occurring on solution branches of
in-phase, anti-phase and phase-locked solutions, respectively. 
\edit{We also keep the corresponding colour-code convention for branches of solutions
and solutions profiles of in-phase, anti-phase and phase-locked type.}

\begin{figure}
  \centering
  \includegraphics[width=\textwidth]{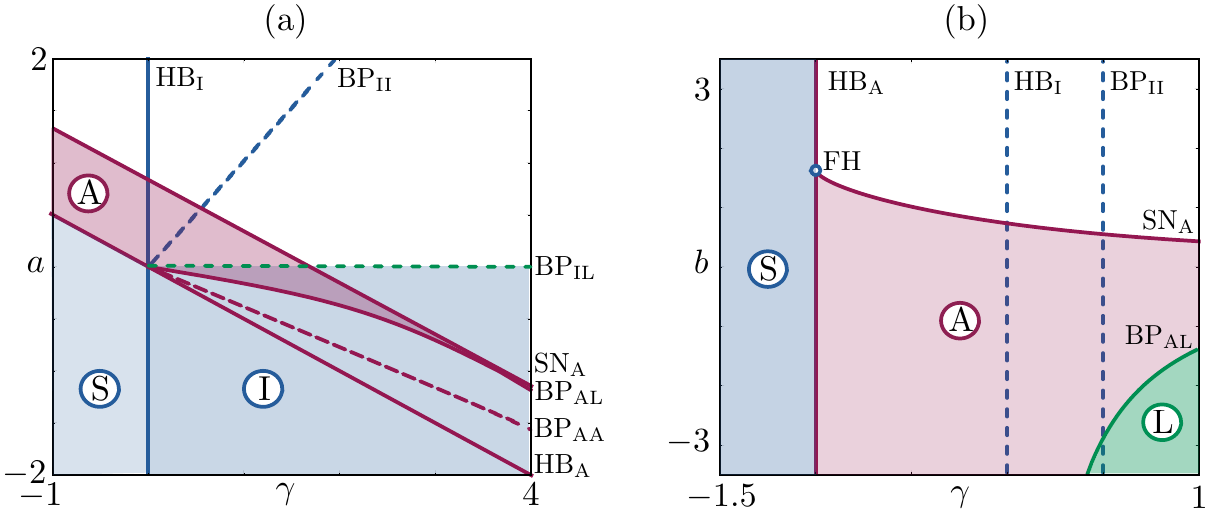}
  \caption{Two parameter continuations of bifurcations occurring in
  Figure~\ref{fig:gammaContinuationCompound}. Panel (a): we fix $b=0.5$ and continue in the
  $(\gamma,a)$-plane the bifurcations of the top panels of
  Figure~\ref{fig:gammaContinuationCompound}; shaded areas represent regions of
  stability for steady states (S), anti-phase (A) and in-phase (I) periodic
  solutions. Panel (b): we fix $a=0.5$ and continue in the $(\gamma,b)$-plane the
  bifurcations in the right panels of Figure~\ref{fig:gammaContinuationCompound}; stable
  phase-locked solutions are indicated by (L).}
  \label{fig:gamma_a_gamma_b_continuation}
\end{figure}

In-phase and anti-phase coordination regimes are born via Hopf bifurcations
($\textrm{HB}_\textrm{I}$, $\textrm{HB}_\textrm{A}$) of the trivial steady state.  In the first
quadrant (where the coupling strength parameters are both positive), anti-phase coordination is the only stable state: a branch of unstable
in-phase solutions is born at $\textrm{HB}_\textrm{I}$ and bifurcates at a
symmetry-breaking bifurcation, $\textrm{BP}_\textrm{II}$, giving rise to a secondary
branch of unstable in-phase solutions where the oscillation amplitudes for agent $1$
and $2$ differ. In the third quadrant (where the coupling strength parameters are both negative) the scenario is specular: in-phase oscillations
are now stable, while anti-phase solutions are unstable and bifurcate at $\textrm{BP}_\textrm{AA}$.
In the first and third quadrants of the $(a,b)$-plane ($a=0.5$, $b=0.5$ and $a=-0.5$,
$b=-0.5$, respectively) there are no branches of phase-locked solutions.

Phase-locked coordination regimes arise in the second and fourth quadrants of the
$(a,b)$-plane, at symmetry-breaking bifurcations of anti-phase solutions
($\textrm{BP}_\textrm{AL}$). Phase-locked solutions are found to be always stable (unstable) in the
fourth (second) quadrant. We note that, in these quadrants, the coupling
nonlinearities $I_{12}$ and $I_{21}$, as functions of $x_1 - x_2$, attain
both negative and positive values as opposed to what happens in the first and third
quadrants, where such functions are strictly positive and strictly negative,
respectively. Stable phase-locked solutions, spanning relative phases in the range of
$(0^{\circ}, 180^{\circ})$, exist for $a>0$ and $b<0$. 
\edit{Such parameter settings could be used to model experiments displaying
coordination regimes different from the canonical in-phase and anti-phase
ones (see the Discussion section for more details).}

\begin{figure}
  \centering
  \includegraphics{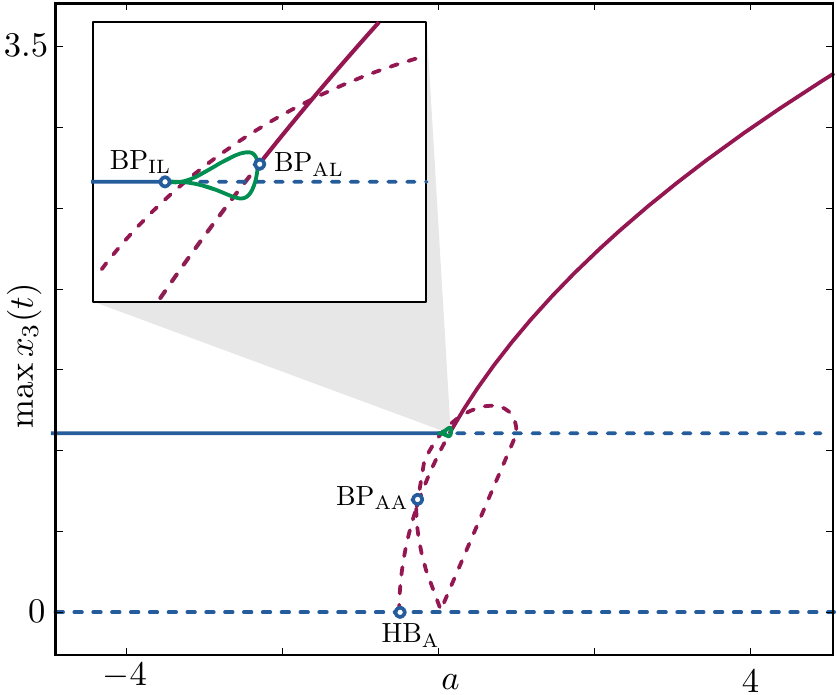}
  \caption{Continuation in the parameter $a$, for $b=-0.5$ and $\gamma=1$, $\omega=2$,
  $\alpha=\beta=1$. The branches show that, with suitable combination of the parameters, it is
  possible to have stable in-phase, anti-phase, and phase-locked oscillations by
  varying $a$. Solid lines represent stable and dashed lines unstable states of \eqref{eqn4}.}
  \label{fig:a_continuation_b_m05}
\end{figure}

\begin{figure}
  \centering
  \includegraphics{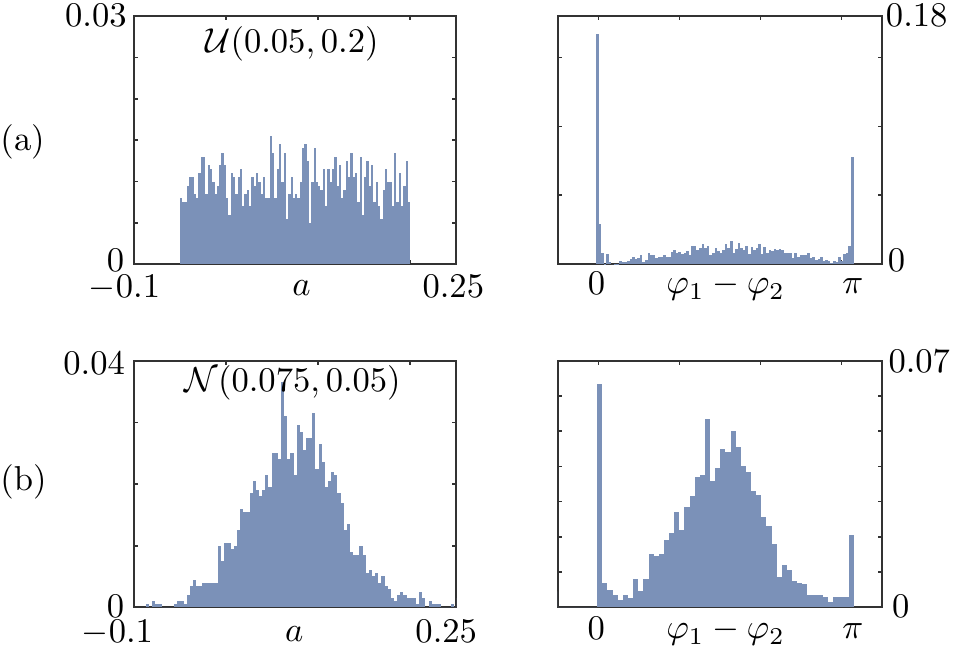}
  \caption{\edit{Distribution of the phase difference between $x_1(t)$ and
  $x_2(t)$ (left panels) obtained when the control parameter $a$ is randomly
  distributed (right panels) with values close to the inset of
  Fig.~\ref{fig:a_continuation_b_m05}. (a): the parameter $a$ is sampled from a
  uniform distribution (left) and 2000 independent simulations are performed; the
  histogram for the relative phase $\varphi_1 - \varphi_2$ is bimodal and sharply
  peaked around $0^\circ$ and
  $180^\circ$, with a small but non-zero probability of finding an intermediate phase
  lag. (b): the experiment is repeated with a normal distribution, which causes a
  third peak to develop around $90^\circ$ in the distribution for the phase lag; the
  latter peak is inherited from the distribution of the control parameter $a$.}}
  \label{fig:uncertaintyQuantificationA}
\end{figure}

\edit{
In Fig.~\ref{fig:gamma_a_gamma_b_continuation} we summarise the behaviour of the
representative examples reported above, for selected values of $a$ and $b$, by
continuing in $(\gamma,a)$ and $(\gamma,b)$ all bifurcation points found in
Fig.~\ref{fig:gammaContinuationCompound}.
}
The two-parameter continuations are
performed so as to show how the solution landscape changes as we pass from the
first to the second quadrant (continuation in $(\gamma,a)$-plane) and from the
second to the third quadrant (continuation in $(\gamma,b)$-plane).
In these two-parameter bifurcation diagrams, we highlight areas where stationary
and oscillatory solutions are stable. In the $(\gamma,a)$-plane, the organising
centre is at $\gamma=0$, $a=0$: at this point, the eigenvalues of the linearised
Jacobian at the trivial state $(x_1,x_2,x_3,x_4)=(0,0,0,0)$ are purely
imaginary, equal to $\pm 2 i$, each with multiplicity $2$, corresponding to
eigenvalues $(0,\mp i/2,0,1)$ and $(\mp i/2,0,1,0)$.
\edit{
For low positive values of the damping $\gamma$, the system supports stable
in-phase solutions (for negative values of $a$) and stable anti-phase solutions (in
a wedge delimited by the locus of $\textrm{SN}_\textrm{A}$ and
$\textrm{BP}_\textrm{AA}$). In a sizeable} region of parameter space, stable in-phase and
anti-phase solutions coexist (see intersection between magenta- and blue-shaded
areas). We note that the original set of parameter values based on \cite{Haken1985} could be found in this region.

In the $(\gamma,b)$-plane, the organising centre is a fold-Hopf bifurcation
around $\gamma=-1$, $b\approx 1.625$ (FH in
Figure~\ref{fig:gamma_a_gamma_b_continuation}) where the locus of saddle-nodes
of the anti-phase solutions, $\textrm{SN}_\textrm{A}$, collides with the locus
of Hopf bifurcations $\textrm{HB}_\textrm{A}$. In this region of parameter
space, phase-locked solutions are found for sufficiently high damping and
sufficiently negative values of $b$. It should be noted, however, that
phase-locked and anti-phase oscillations do not coexist for the default choice of intrinsic dynamics parameter values. 
As it can also be verified analytically, the locus of bifurcations $\textrm{HB}_\textrm{A}$ and
$\textrm{HB}_\textrm{B}$ of the stationary steady state do not depend on $b$.
It should be also noted that, for suitable combination of the parameters, it is possible
to \edit{visit} stable in-phase, anti-phase, and phase-locked oscillations by varying $a$.
An example of continuation in $a$ for $b=-0.5$ and $\gamma=1$, $\omega=2$, $a=b=1$ is
depicted in Figure~\ref{fig:a_continuation_b_m05}. 
\edit{
It can be clearly seen in the
inset of Fig.~\ref{fig:a_continuation_b_m05}
that, as $a$ increases, the stable in-phase coordination regime (characterised by
relative phase $0^{\circ}$) loses stability; then, a phased-locked coordination regime
emerges (ranging over relative phases in $(0^{\circ}, 180^{\circ})$ in a continuous
fashion) and eventually a stable anti-phase coordination regime (characterised by
relative phase $180^{\circ}$) is established.}
\edit{The bifurcation diagram implies that, in an experimental setup where
$a$ were to be assigned randomly, we would observe trajectories with relative phases
distributed in the interval $(0^{\circ}, 180^{\circ})$, and with peaks at $0^\circ$
and $180^\circ$. To verify this prediction, we performed an uncertainty
quantification study, in which $a$ is assigned randomly, near the inset of
Fig.~\ref{fig:a_continuation_b_m05}, and histograms of relative phases between
$x_1(t)$ and $x_2(t)$ are computed a posteriori. In
Fig.~\ref{fig:uncertaintyQuantificationA}(a) we perform $2000$ independent
simulations, where $a$ is sampled from the uniform distribution between 0.05 and 0.2
and plot the resulting phase lag histogram: the distribution for the phase difference
$\varphi_1 - \varphi_2$ (right) is bimodal and sharply peaked around $0^\circ$ and
$180^\circ$, as
expected, with a
small but non-zero probability of finding an intermediate phase lag. The likelihood
that experiments display such intermediate relative phases is deeply affected by the
distribution of $a$: if we pass from a uniform to a normal distribution for $a$,
(Fig.~\ref{fig:uncertaintyQuantificationA}(b)), the
resulting phase lag distribution develops also a peak around $90^\circ$, which is
inherited from the parameter distribution.}

\begin{figure}
  \centering
  \includegraphics[width=\textwidth]{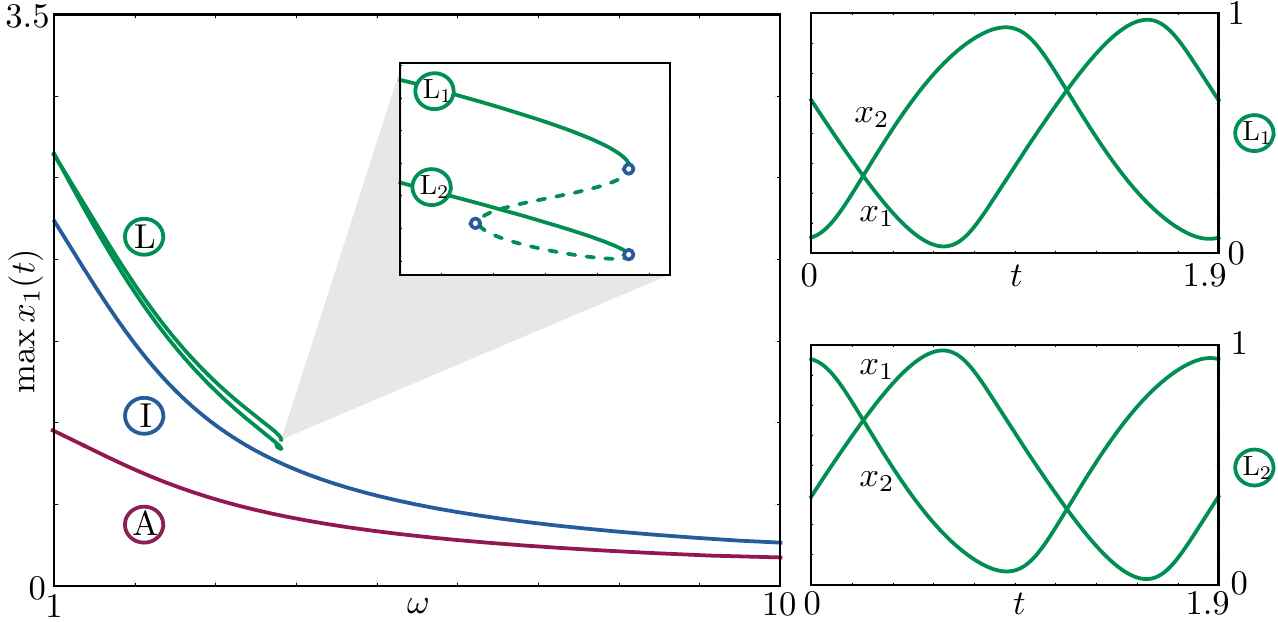}
  \caption{Continuation of in-phase, anti-phase and phase-locked
  solutions in the frequency $\omega$. Solutions behave similarly to the single
  HKB oscillator case
  (Figure~\ref{fig:singleOscillator}b), except they have various phase
  behaviours. The branch of phase-locked solutions undergoes a series of
  saddle-nodes bifurcation, giving rise to stable solutions in which the
  phase difference is reversed (see solutions profiles $\mathrm{L}_{1,2}$). We
  note that the branches in this figure do not coexist, as they are found in
  different regions of parameter space: $\alpha=\beta=1$ and 
  $a=-0.5$, $b=-0.5$, $\gamma=5$ (in-phase), 
  $a= 0.5$, $b=-0.5$, $\gamma=1.2$ (anti-phase) and $a= 0.5$, $b=-0.5$,
  $\gamma=6.2$ (phase-locked).}
  \label{fig:omegaContinuation}
\end{figure}

It is interesting to study the behaviour of various periodic solutions as the
common eigenfrequency of the oscillators, $\omega$, varies. We selected stable in-phase,
anti-phase and phase-locked solutions and continued them in $\omega$
(Figure~\ref{fig:omegaContinuation}). We found that such solutions behave
essentially as in the single oscillator case (Figure~\ref{fig:singleOscillator}b): low
frequencies elicit large-amplitude oscillations with abrupt time transitions,
whereas large frequencies induce smoother small-amplitude oscillations. In this
case the changes in the oscillation patterns occur to both agents, with various
phase differences. The branch of phase-locked solutions undergoes a series of
saddle-node bifurcations, giving rise to stable solutions in which the phase
difference is reversed (see solutions profiles $\mathrm{L}_{1,2}$ in
Figure~\ref{fig:omegaContinuation}).

\begin{figure}
  \centering
  \includegraphics{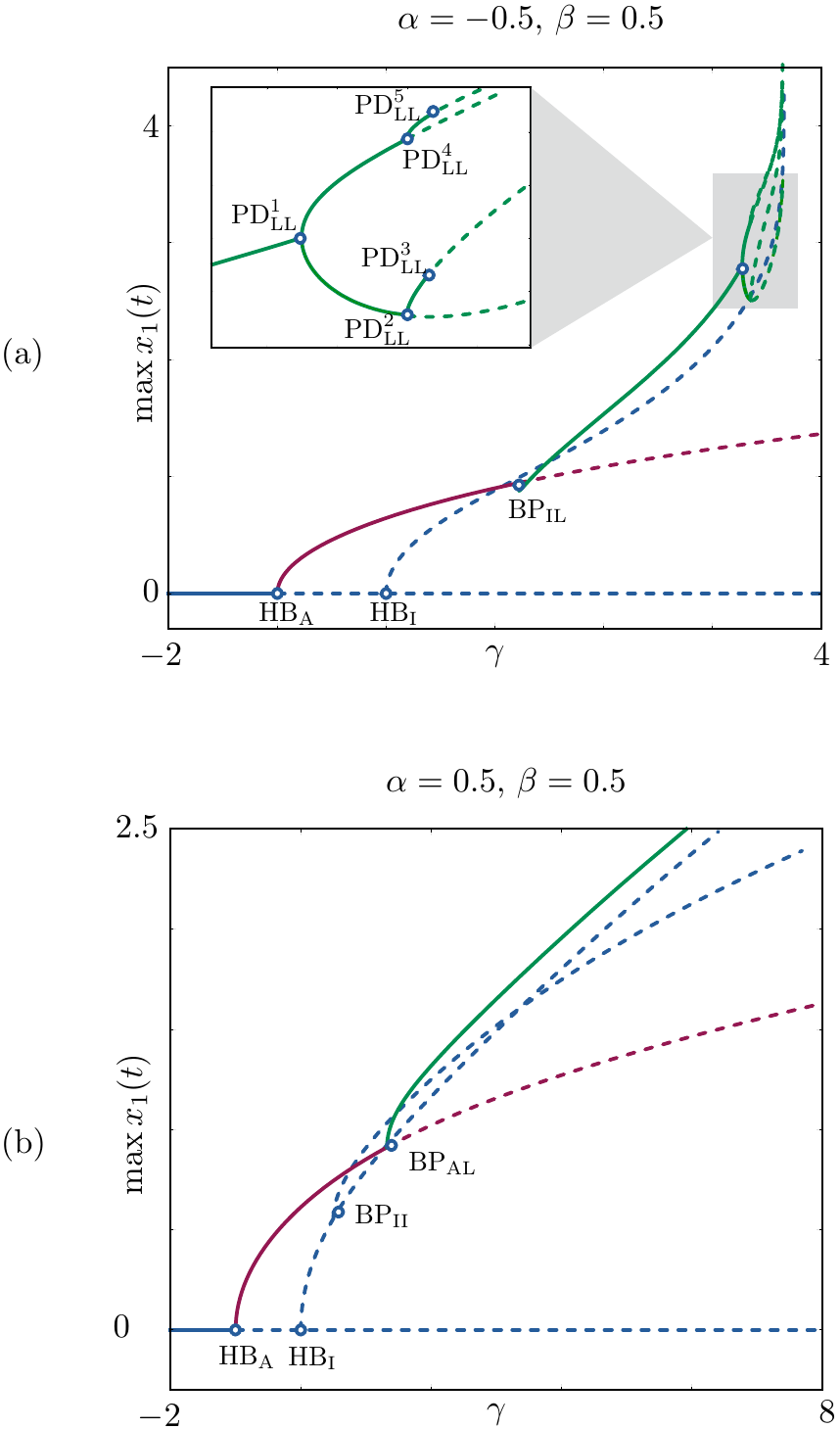}
  \caption{Bifurcation diagram in $\gamma$ for $\omega=2$, $a = 0.5$, $b=-0.5$
  and various values of $\alpha$ and $\beta$. Solid lines represent stable and dashed lines unstable states of \eqref{eqn4}.}
  \label{fig:gammaContinuationCompound_alpha_beta}
\end{figure}

\begin{figure}
  \centering
  \includegraphics{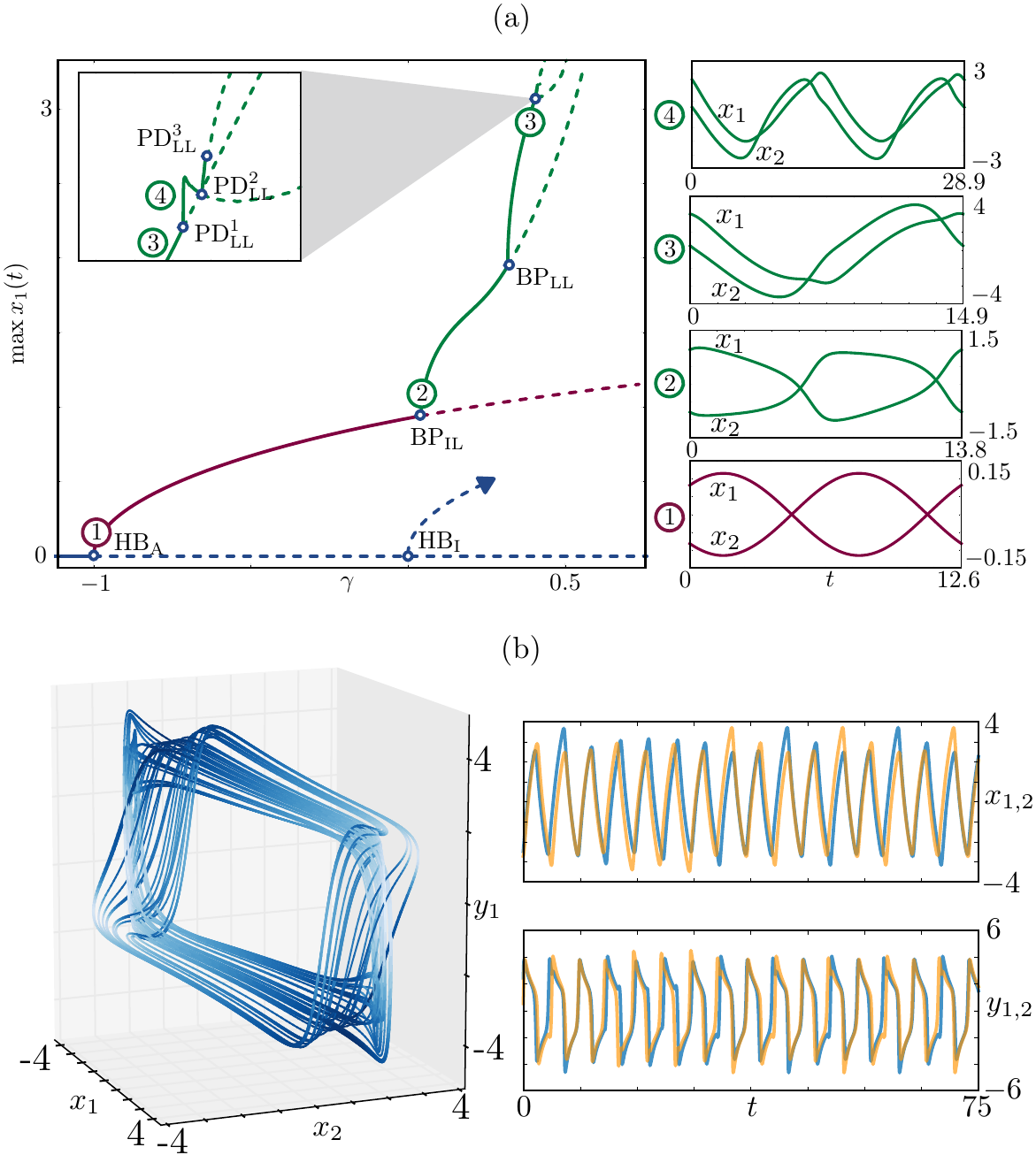}
  \caption{Period doubling cascade. (a): branch with Period Doubling cascade and
  stable non-trivial periodic solutions. We plot one period of several stable
  solutions along the branch, whereas we omit the unstable branch emanating from
  \edit{$\textrm{HB}_\textrm{I}$}. Solutions feature increasing solution periods, $T_1
  \approx 12.6$, $T_2 \approx 13.8$, $T_3 \approx 14.9$, $T_4 \approx 28.9$,
  corresponding to $\Omega_1 \approx 0.500$, $\Omega_2 \approx 0.455$, $\Omega_3
  \approx 0.422$, $\Omega_4 \approx 0.217 $, respectively. Parameters
  $\omega=0.5$, $a = 0.5$, $b=-0.5$ $\alpha=0.5$ and
  $\beta=-0.05$. (b): attractor found for $\gamma=3.42$; the coordination regime
  shows erratic phase changes, during which $x_1$ and $x_2$ alternate in the
  leading position. This regime involves fast velocity switches, as evidenced by
  the time traces of $y_1$ and $y_2$. Solid lines represent stable and dashed lines unstable states of \eqref{eqn4}.}
  \label{fig:PDExample}
\end{figure}

Finally we investigate the impact of intrinsic oscillator dynamics on the collective behaviour of the HKB model by performing bifurcation analysis in the intrinsic dynamics parameters
$\alpha$ and $\beta$. Instead of presenting two-parameter bifurcation diagrams for
different cases, we report here only notable examples of our computations (see
Figures~\ref{fig:gammaContinuationCompound_alpha_beta} and
\ref{fig:PDExample}(a)). The bifurcation structures found in these cases have
common traits with the ones discussed above for the coupling strengths parameters $a$ and $b$, that
is, the trivial steady state undergoes Hopf bifurcations to anti-phase and
in-phase periodic states, and various symmetry-breaking bifurcations give rise
to phase-locked solutions. Interestingly, when varying $\alpha$ and $\beta$ we
could find period doubling cascades, which are found robustly when $\alpha$ and
$\beta$ have opposite signs, as evidenced in 
Figure~\ref{fig:gammaContinuationCompound_alpha_beta}(a), where $\alpha=-0.5$,
$\beta= 0.5$, and Figure~\ref{fig:PDExample}(a), where $\alpha=0.5$, $\beta= -0.05$.
Representative stable solutions on the period-doubling cascade are also shown in
Figure~\ref{fig:PDExample}(a). 

Using direct numerical simulations we explored the system behaviour close to the
period-doubling cascade, finding chaotic regimes (see Figure~\ref{fig:PDExample}(b)) in
which the solution remains bounded and features sudden erratic phase transitions, during which
the agents alternate as leaders and followers. In this regime, the velocities
$y_1$ and $y_2$ undergo fast switches. 
\edit{The existence of such complex solutions is perhaps not surprising from a dynamical
systems viewpoint; however, the behaviour described above has not been reported
nor investigated previously, and can be used to model experiments where the
movement coordination is irregular in nature.
Last but not least, knowledge about the existence of such solutions is critical when
designing virtual player interaction environments
\cite{Kelso2009,ZhaiSMC2014,ZhaiCDC2014,Zhai2016} and/or planning human dynamic clamp
experiments based on the HKB model \cite{Dumas2014}.}

\subsubsection{Bistability and hysteresis}

\begin{figure}
  \centering
  \includegraphics[width=\textwidth]{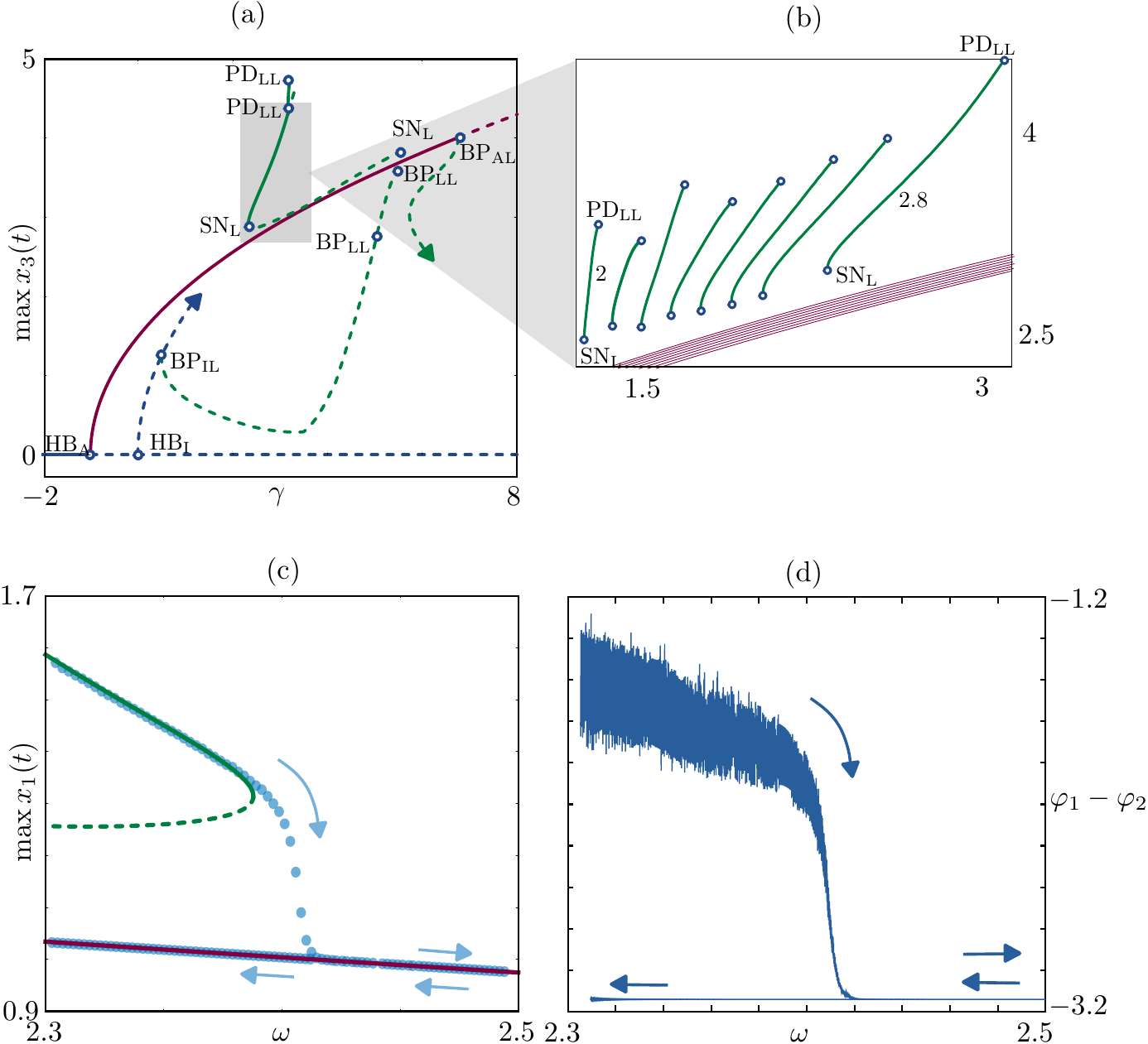}
  \caption{Bistability and hysteresis between anti-phase and phase-locked
  solutions. (a): the in-phase branch (blue) undergoes a
  \edit{symmetry-breaking bifurcation ($\textrm{BP}_\textrm{IL}$) and the
  resulting unstable phase-locked branch, featuring two further symmetry-breaking
  bifurcations ($\textrm{BP}_\textrm{IL}$), restabilises at a saddle node
  bifurcation, before a period-doubling cascade takes place.}
  A stable portion of the
  phase-locked branch (solid green line between $\textrm{SN}_\textrm{L}$ and
  $\textrm{PD}_\textrm{LL}$), coexists with the anti-phase branch originating at
  $\textrm{HB}_\textrm{I}$ (solid red branch). Parameters: $a=0.5$, $b=0.5$,
  $\omega=3$, $\alpha=-1.7$, $\beta=0.5$. (b): we repeat the experiment for
  $\omega \in [2,2.8]$ and plot stable branches to highlight the bistability
  region. (c): $\omega$ is varied by continuation and by quasi-static sweeps in
  direct numerical simulations (blue dots), for $\gamma=1.7$; the time
  simulation follows the phase-locked branch up to the saddle node at $\omega
  \approx 2.4$, where
  an abrupt and hysteretic transition to an anti-phase solution is observed.
  (d): phase lag during numerical simulation in (c).} 
  \label{fig:bistabilityHysteresisCompound}
\end{figure}

In this section we explore further the dependence of the HKB model dynamics on the intrinsic properties of the coupled oscillators.
In suitable regions of parameter space we find coexisting stable periodic states
characterised by different relative phases or phase lags. In
Figure~\ref{fig:bistabilityHysteresisCompound}(a) we run a
continuation similar to the ones presented above, but we set
$\alpha=-1.7$. The branches of this bifurcation diagram are qualitatively
similar to the ones of the previous sections; however the in-phase periodic
branch originating at the Hopf bifurcation $\textrm{HB}_\textrm{I}$ undergoes a
symmetry-breaking bifurcation ($\textrm{BP}_\textrm{IL})$. Such branch is
initially unstable, undergoes 2 other symmetry-breaking bifurcations,
re-stabilises at a saddle node bifurcation and then features a period doubling
cascade. The stable portion of this branch (solid green branch between
$\textrm{SN}_\textrm{L}$ and $\textrm{PD}_\textrm{LL}$) coexists with a branch
of stable anti-phase solutions originating from the trivial state at
$\textrm{HB}_\textrm{A}$ (red branch).

\edit{
This bifurcation structure opens up the possibility of observing abrupt
relative phase transitions between phase-locked (at any relative phase between
$0^{\circ}$ and $180^{\circ}$) and anti-phase (at relative phase
equal to $180^{\circ}$) coordination regimes as a function of the eigenfrequency $\omega$.}
We find that bistability is observed in a significant region of parameter space: in the
inset of Figure~\ref{fig:bistabilityHysteresisCompound}(b) we report overlapping
stable portions of phase-locked and anti-phase branches as we vary the
eigenfrequencies $\omega$. As $\omega$ is varied, the anti-phase branch (red) changes only
slightly, while the stable phase-locked branch moves to the left and expands.
Our analysis predicts coexistence in the region $(\gamma,\omega) \in
[1.2,3.2]\times[2,2.8]$ (which was found robustly in other parameter regions, not
shown). It is important to note that this phase transition is qualitatively different
from the transition addressed by the original HKB model \cite{Haken1985}, where an
increase in frequency leads to transition from anti-phase to in-phase coordination.
In the parameter regime described above, an increase in frequency leads to transition
from phase-locked to anti-phase coordination behaviour.

\begin{figure}[H]
  \centering
  \includegraphics[width=\textwidth]{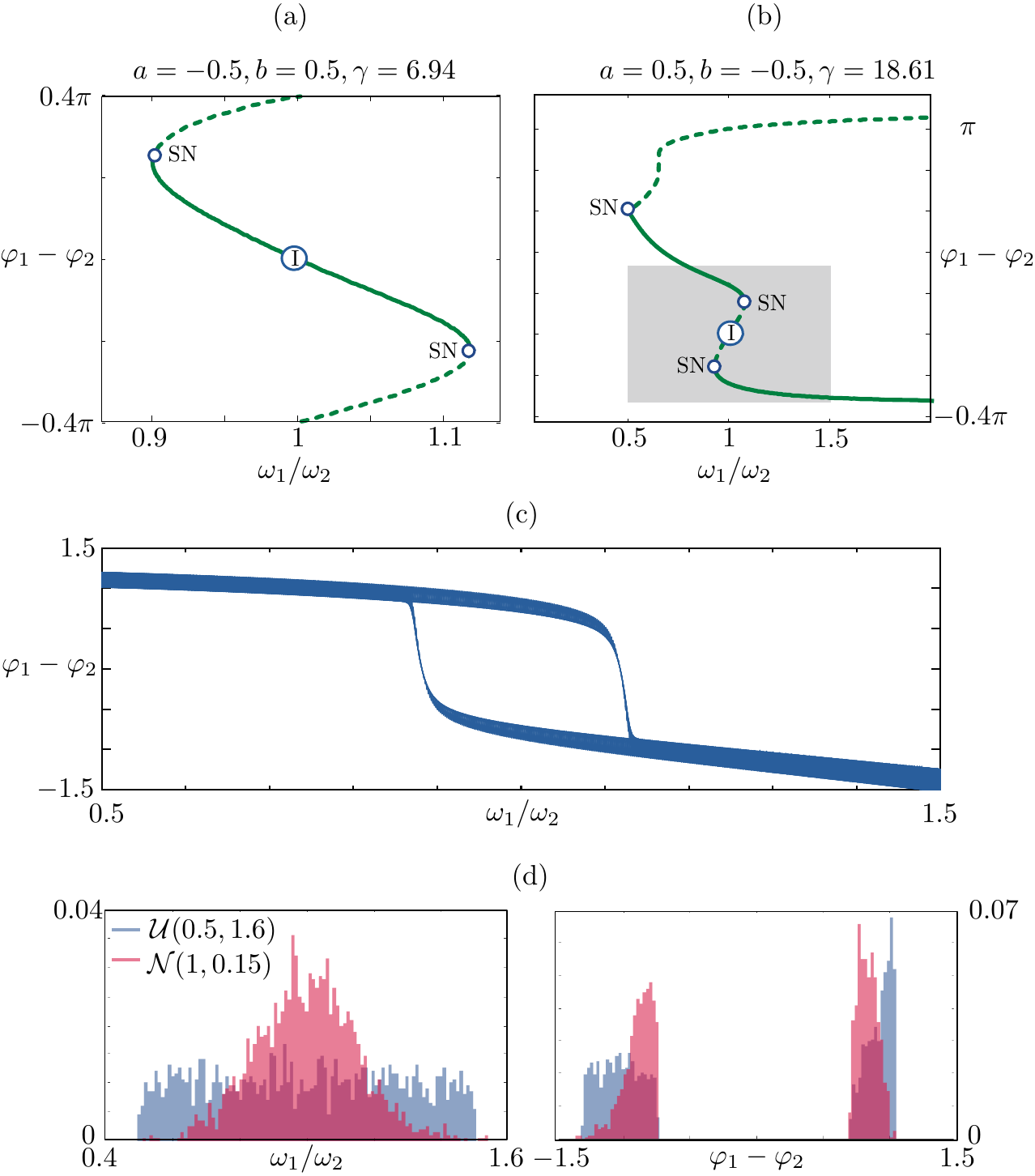}
  \caption{Phase difference between periodic solutions $x_1(t)$ and $x_2(t)$ as
  a function of the ratio $\omega_1/\omega_2$. (a): A stable solution on the
  in-phase branch in the fourth quadrant of
  Figure~\ref{fig:gammaContinuationCompound} is continued in $\omega_1/\omega_2$.
  (b): The continuation is repeated starting from a solution on the phase-locked
  branch in the second quadrant of Figure~\ref{fig:gammaContinuationCompound}. In
  both cases, heterogeneity in the eigenfrequencies impacts the phase lag of the
  solution.
  (c): when the ratio $\omega_1/\omega_2$ is modulated with a slowly-varying
  sinusoidal function, the actors alternate in the leading position with
  hysteretic cycles, which follow the branches in the inset of (b) and jump at
  the corresponding saddle-node bifurcations. \edit{(d): we perform an
  experiment similar to the one in Fig.~\ref{fig:uncertaintyQuantificationA}; when
  the parameter $\omega_1/\omega_2$ is drawn randomly near the shaded area in (b)
  from a uniform (blue) or a normal (red) distribution, the resulting phase lag
  distribution is bimodal, with peaks at $\pm 57.30^\circ$, as predicted by the
  bifurcation diagram in (b) and by the parameter sweep in (c).}} 
  \label{fig:phaseLagVersusOmegaRatio}
\end{figure}

\newpage
To illustrate the dynamical switch between solution types, we perform time-stepping
simulations in which the eigenfrequency $\omega$ is varied quasi-statically and
compare with the bifurcation analysis. In
Figure~\ref{fig:bistabilityHysteresisCompound}(c) we
continued an anti-phase (red) and phase-locked (green) solution for
$\gamma=1.7$, $\omega=2.3$ in the parameter $\omega$; the phase-locked branch
destabilises at a saddle-node bifurcation, whereas the anti-phase branch remains
stable for $\omega \in [2.3,2.5]$. We then initialised a time simulation on the
phase-locked branch (blue dots in
Figure~\ref{fig:bistabilityHysteresisCompound}(c)) and changed
$\omega$ in slow increments (from $\omega = 2.3$ up to $\omega =2.5$) followed
by small decrements (from $\omega = 2.5$ down to $\omega =2.3$). The time
simulation shows an abrupt and hysteretic change in the solution type. This could be 
further appreciated in Figure~\ref{fig:bistabilityHysteresisCompound}(d) where
we plot the time simulation using the phase lag $\varphi_1 - \varphi_2$. The
$x_2$ is delayed with respect to $x_1$, with an initial phase lag
$\varphi_1-\varphi_2 \approx 90^{\circ}$; when $\omega \approx 2.4$, we
observe a transition to an orbit with $\varphi_1 - \varphi_2 \approx 180^{\circ}$
(anti-phase solution).

\subsubsection{Effect of heterogeneity in eigenfrequencies on the coordination
regimes}

In the computations shown so far, the two oscillators possess a common eigenfrequency $\omega$. 
In order to study the effect of heterogeneities on
coordination, we introduce two parameters $\omega_1$, $\omega_2$, then fix $\omega_1$
to the nominal value $\omega_1=2$ and use the ratio $\omega_1/\omega_2$ as a
continuation parameter. 
\edit{
The difference in eigenfrequencies introduces a heterogeneity
in the system and has the potential to turn in-phase solutions into phase-locked
solutions and vice-versa. In order to illustrate this idea we performed bifurcation
analysis in theparameter $\omega_1/\omega_2$ investigating the in-phase solutions
which exist for parameter values within the range of those used in the original HKB
model \cite{Haken1985}, as well the stable phase-locked solutions which we reported above
for $a>0$ and $b<0$ (see Figure~\ref{fig:gammaContinuationCompound}).
}
In Figure~\ref{fig:phaseLagVersusOmegaRatio}, we
initialise the continuation with an in-phase and a phase-locked periodic
solution. We plot the bifurcation diagram in terms of the phase lag (measured in
radians), by computing the approximate phases times $\varphi_i = t_i/T$, for
$i=1,2$, where $t_i$ is the time at which the orbit $x_i(t)$ attains its maximum
and $T$ is the solution period. Figure~\ref{fig:phaseLagVersusOmegaRatio}(a) depicts a stable,
initially in-phase, solution at $\omega_1/\omega_2=1$ that turns into a phase-locked solution as
$\omega_1/\omega_2$ is increased/decreased, losing stability at saddle-nodes bifurcation. In
Figure~\ref{fig:phaseLagVersusOmegaRatio}(b) we show how the phase lag is reduced when
the frequency ratio is varied and an in-phase (albeit unstable) solution is
eventually attained, before a new phase-locked solution arises. 

The bifurcation structure in Figure~\ref{fig:phaseLagVersusOmegaRatio}(b)
implies that hysteresis between phase-locked solutions with opposite phase lags (relative phases)
is possible in the model. To illustrate this we perform
time-stepping simulations in which the ratio is varied quasi-statically as
$\omega_1(t) = \omega_2(t)[1 + \sin(0.005 t)]$ and plot the results in Figure~\ref{fig:phaseLagVersusOmegaRatio}(c). 
The two oscillators swap in the leader and follower role, following the branches of Figure~\ref{fig:phaseLagVersusOmegaRatio}(b) and switching roles at the
corresponding saddle-nodes bifurcations. \edit{This numerical experiment could be interpreted in the light of the joint-improvisation scenario in the "mirror game", 
a recently proposed paradigm for studying the dynamics of two people improvising motion together \cite{Noy2011}. In particular, as the participants are  asked
to imitate each other and create synchronised and interesting motions they would be naturally trying to adjust their movement velocities and thus eigenfrequencies
to each other. This would lead to variation in the ratio of their eigenfrequencies and respectively exchange of leader and follower roles while playing the game.
Indeed, observations based on our data collected in a "mirror game"
setting \cite{Slowinski2016} indicate that the distribution of relative phase during
a typical joint improvisation sessions are bi-modal 
pointing to \edit{possible} hysteretic dynamics.
As we see in Fig.~\ref{fig:phaseLagVersusOmegaRatio}(d), the bimodal
distribution emerges also in the case of randomly assigned frequency ratio
$\omega_1/\omega_2$: when the value of this parameter is drawn randomly (close to the
hysteretic region) from a uniform or a normal distribution, the resulting phase lag
distribution is bimodal, with peaks at $\pm 57.30^\circ$, as predicted by the
bifurcation diagram in Fig.~\ref{fig:phaseLagVersusOmegaRatio}(b) and by the
parameter sweep in Fig.~\ref{fig:phaseLagVersusOmegaRatio}(c).}

\section{Discussion}

In this paper we have systematically investigated the dynamics of the HKB model in the state space spanned by the position and velocities of the coupled oscillators. \edit {Furthermore we go beyond the weakly coupled regime and consider the coupling strength parameters as generic.}
We show that stable periodic solutions in the single HKB oscillator model are born
via a Hopf bifurcation as the damping parameter $\gamma$ becomes positive.
Furthermore we reveal that, under certain intrinsic oscillator properties the
periodic solutions of the single HKB model could disappear via a heteroclinic cycle,
associated with rapid increase in the magnitude of the state
variables.
Although, such behaviour cannot be observed in a physical system it can have significant consequences for the design and development of the virtual players.
Bifurcation analysis of the full four dimensional HKB model reveals a variety of
different coordination dynamics. Attractors at a constant relative phase of
$0^{\circ}$ (in-phase) and/or $180^{\circ}$ (anti-phase) are born via Hopf
bifurcations detected in the damping parameter $\gamma$. 
We find symmetrical attractors of phase-locked solutions at intermediate values of
relative phase (between $0^{\circ}$ and $180^{\circ}$), which increase or decrease
gradually as $\gamma$ is increased. \edit{We demonstrate that the phase-locked solutions
are born in a symmetry-breaking bifurcation of periodic orbit in which the anti-phase
periodic attractor loses stability as the damping parameter $\gamma$ is varied.}
Changing the sign of the coupling strengths has the effect of shifting the attractors
by $180^{\circ}$ thus changing the phase that remains stable at high frequencies,
from $0^{\circ}$ to $180^{\circ}$ or vice-versa.
We also show that change in the intrinsic oscillators' properties (i.e. when varying
the parameter $\alpha$) can lead to complex dynamics mediated via a period-doubling
cascade. Furthermore different intrinsic dynamics can also bring about a variety of
bi-stability modes, which are different that the type of bi-stability described in
the original HKB model study \cite{Haken1985}.
Finally we consider a case of a heterogeneity in the system by introducing difference
in the eigenfrequency of the coupled oscillators. We demonstrate how this results in
bi-stability and hysteresis. \edit{Our uncertainty quantification simulations presented in Fig.~\ref{fig:phaseLagVersusOmegaRatio}(d) confirm that in the case of heterogeneous oscillators hysteresis loops and phase-locked coordination modes should be expected in experiments, as suggested in \cite{Bardy2002}. What is more, existence of such hysteresis loop provides an excellent opportunity for a quantitative experimental validation of the HKB model using two heterogeneous coupled oscillators e.g. by putting weights on body parts as suggested in \cite{Bardy2002} or by using heterogeneous pendula as in \cite{Schmidt2008}}.

In a large number of multi-stable examples observed experimentally the patterns of stability change under
different conditions. Bimanual finger coordination is bi-stable at low frequencies
but above a critical frequency the anti-phase pattern is no longer sustainable
\cite{Haken1985}. Similarly postural sway is bistable at low frequencies ($20^{\circ}$ and $180^{\circ}$) but
the phase-locked ($20^{\circ}$) mode looses stability at high frequencies or when other behaviours,
such as reaching, are incorporated in the task \cite{Bardy2002}. These transitions
between stable states, and particularly the loss of stability of the anti-phase
mode at high frequencies, appear to be a fundamental feature of human coordination
\cite{Kelso1997}. The hypothesis that these real-world patterns and transitions between them
are emergent phenomena due to a self organised dynamical system are substantiated
by experimental results such as critical fluctuations, critical slowing
down and hysteresis between modes \cite{Bardy2002,Fuchs2008},\cite{Schoener1986}.
In this paper we make the first step towards identifying parameter regimes and dynamics that would allow
to model a variety of different experimental observations using the same modelling framework.

Many recent experimental studies of human movement coordination 
\cite{Issartel2007,Bourbousson2010,Zanone1997,Esteves2012,Duarte2012,Davids2013} have reported
persistent movement coordination dynamics other than the well known in-phase and
anti-phase synchronisation behaviour that have inspired the development of the HKB
model \cite{Haken1985}. 
Despite the large number of behaviours whose dynamics are well represented
by the theoretically predicted in-phase and anti-phase stability, there are several
counter examples where~\edit{human body movements show evidence of stability at different or additional
intermediate phases}. Examples of real-world systems with stabilities at other
relative phases include: the human postural system (stability at $20^{\circ}$)
\cite{Bardy2002}, amble to walk
gait in quadrupeds (stability at $90^{\circ}$) \cite{Collins1993}, the bipedal skipping gait \cite{Minetti1998}, coordination
tendencies of successful defences ($30^{\circ}$) and unsuccessful defences ($90^{\circ}$) in
soccer \cite{Duarte2012}, squash ($135^{\circ}$) \cite{McGarry2006}, and butterfly stroke swimming ($90^{\circ}$) \cite{Ehrlacher2003} as well as variety of relative phase distributions in other team sports \cite{Davids2013} . There is evidence that other phases can be stable simultaneously with $0^{\circ}$ and $180^{\circ}$. These multi-stable dynamics can exist naturally or be learnt \cite{Zanone1997}. Our results about the existence of stable phased-locked periodic solutions in the HKB model that span all possible relative phases between $0^{\circ}$ and $180^{\circ}$ could be related to some of the above mentioned experimental observations.
In particular, analysis of the data collected from interactions between player dyads allowed for a description of the space-time dynamics of basketball match-play. In the longitudinal direction, a strong attraction to in-phase was reported for all possible dyads but not so for the lateral direction. Instead, attractions to in-phase or anti-phase were observed among most dyads with the player vs. player dyads tending on balance to demonstrate less pronounced attractions or repulsions to certain relative-phases than the player-opponent dyads \cite{Bourbousson2010,Esteves2012}. Interpersonal coordination tendencies of 1-vs-1 sub-phases were investigated in \cite{Duarte2012}. 
The experimental results presented in Fig. 2 in \cite{Duarte2012} could be, for
example, qualitatively accounted for by the type of coordination stability dependence
on the coupling parameter $a$ found in the HKB model. \edit{Specifically (see
Fig. \ref{fig:a_continuation_b_m05} for $b=-0.5$ and $\gamma=1$, $\omega=2$,
$\alpha=\beta=1$), as the coupling strength between the velocity components of the
two oscillators increases, the stable coordination regime exhibited by the HKB model
undergoes transitions form stable in-phase coordination for $a<0$ through stable
phase-locked coordination spanning all possible relative phases in
$(0^{\circ},180^{\circ}$); as $a$ increases, a stable anti-phase coordination regime
with increasing amplitudes is generated.}

\edit{Last but not least, very recent experiments involving the use of virtual
partner interaction \cite{Kelso2009,Dumas2014,ZhaiSMC2014,ZhaiCDC2014,Zhai2016} have
employed to various degree the HKB model in order to study social interactions and
interpersonal coordination. These studies have used adaptation in the HKB parameter
values in their implementations. Knowledge about how the type and stability of the
possible HKB model solutions depend on the model parameters could greatly facilitate
the design and ensure robustness of such hybrid systems where a human interacts with
a virtual partner whose movements are driven by the HKB model. Furthermore,
comparison of the theoretical predictions and dynamics observed in experiments with a
virtual partner could allow for quantitative, rather than qualitative, validation of
different models of motor coordination. Although deficits of the HKB model are
well-known, see for example discussion in \cite{Beek2002}, our analysis demonstrates
that this model has much richer dynamics than previously considered and showcases
mathematical tools that could be very useful in future studies of human movement
coordination.}

\section*{Acknowledgments}
The authors would like to thank Ed Rooke (University of Bristol, UK) for initial
discussions on the HKB model analysis and Pablo Aguirre (Universidad Tecnica Federico
Santa Maria, Chile) for helpful discussions on the global transition in the single
HKB oscillator.  This work was funded by the European Project AlterEgo FP7 ICT 2.9 -
Cognitive Sciences and Robotics, Grant Number 600610. The research of KT-A was
supported by grants EP/L000296/1 and EP/N014391/1 of the Engineering and Physical
Sciences Research Council (EPSRC).

\bibliographystyle{plain}
\bibliography{da_ps_bb_kt_arxiv}

\end{document}